\documentclass[12pt]{article}
\font\smallit=cmti10
\font\smalltt=cmtt10
\font\smallrm=cmr9

\usepackage{epsfig}

\newenvironment{packed_enumerate}{
\setlength{\parsep}{0pt}
\setlength{\parskip}{0pt}
\begin{enumerate}
  \setlength{\itemsep}{1pt}
  \setlength{\parsep}{0pt}
  \setlength{\parskip}{0pt}
}{\end{enumerate}}

\begin{document} 

\begin{center}
\vspace*{-20pt} 
\centerline{\smalltt INTEGERS: \smallrm ELECTRONIC JOURNAL OF COMBINATORIAL NUMBER THEORY
\smalltt 7 (2007), \#G01} 
\vskip 20pt

\uppercase{\bf Diagonal Peg Solitaire}
\vskip 20pt
{\bf George I. Bell\footnote{\tt http://www.geocities.com/gibell.geo/pegsolitaire/}}\\
%{\smallit 5040 Ingersoll Pl, Boulder, CO 80303, USA}\\
{\smallit Tech-X Corporation, 5621 Arapahoe Ave Suite A, Boulder, CO 80303, USA}\\
{\tt gibell@comcast.net}\\
\end{center}
\vskip 30pt
\centerline{\smallit Received: 12/19/05, Revised: 11/7/06, Accepted: 12/23/06, Published: 1/24/07}
\vskip 30pt 

\centerline{\bf Abstract}
\noindent
We study the classical game of peg solitaire when diagonal jumps are allowed.
We prove that on many boards, one can begin from a full board
with one peg missing, and finish with one peg anywhere on the board.
We then consider the problem of finding solutions that minimize the
number of moves (where a move is one or more jumps by the same peg),
and find the shortest solution to the ``central game",
which begins and ends at the center.
In some cases we can prove analytically that our solutions are the
shortest possible,
in other cases we apply A* or bidirectional search heuristics.

\pagestyle{myheadings}
\markright{\smalltt INTEGERS: \smallrm ELECTRONIC JOURNAL OF COMBINATORIAL NUMBER THEORY
\smalltt 7 (2007), \#G01\hfill}

\thispagestyle{empty} 
\baselineskip=15pt 
\vskip 30pt

\section*{\normalsize 1. Introduction}
\begin{figure}[htbp]
\centering
\epsfig{file=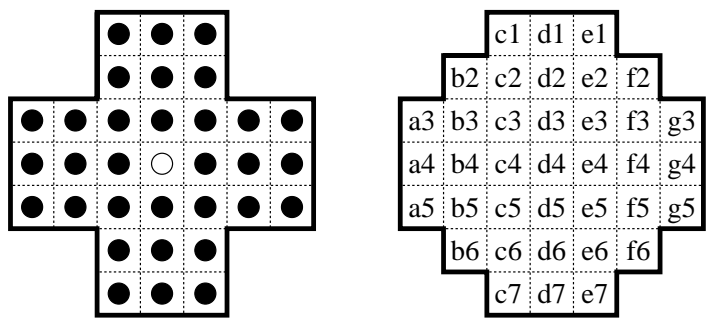}
\caption{(a) The standard 33-hole board. (b) The 37-hole board, with hole coordinates.}
\label{fig1}
\end{figure}

Peg solitaire is a puzzle that has been popular for over 300 years;
it is most commonly played on the 33-hole or 37-hole
boards of Figure~\ref{fig1}.
We refer to a board location as a hole, because on an actual board
there is a hole or depression in which the peg (or marble) sits.
The game begins with pegs in all the holes except one (Figure~\ref{fig1}a).
The player jumps one peg over another into an empty hole,
removing the jumped peg from the board.
The goal is to select a sequence of jumps that finish with one peg.

In the standard version of the game,
only horizontal and vertical jumps are allowed (along rows and columns),
and there can be at most four different jumps into any hole.
Following Beasley \cite{Beasley}, we will refer to this game
as \textbf{4-move} solitaire.
In this paper, we'll explore the version of the game where
diagonal jumps in both directions are also allowed---there can then be
up to eight jumps into a hole and this will be called \textbf{8-move}
solitaire.
An intermediate version, in which diagonal jumps are allowed in only one direction,
is called \textbf{6-move} solitaire, and is equivalent to
solitaire played on a triangular grid \cite[p. 233]{Beasley}.

\begin{figure}[htbp]
\centering
\epsfig{file=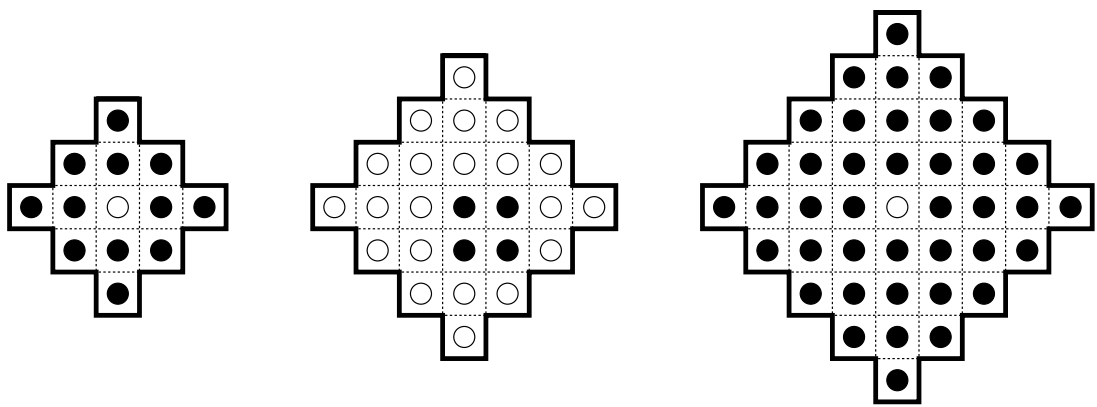}
\caption{Diamond boards 3, 4, and 5 holes on a side with 13, 25, and 41 holes.}
\label{fig2}
\end{figure}

In this paper, we will consider 8-move solitaire on five 
square symmetric boards: the ``standard" 33-hole board
(Figure~\ref{fig1}a), the 37-hole board
(Figure~\ref{fig1}b), plus three diamond-shaped boards
of various sizes (Figure~\ref{fig2}).
The board $Diamond(n)$ has $n$ holes on a side, and a total of
$n^2+(n-1)^2$ holes.
The \textbf{central game} is the problem which begins with a full-board
with one peg missing at the center, and finishes with one peg in the center.
In 4-move solitaire, the central game is not solvable on any of these
boards except for the standard 33-hole board.
However, as we shall see, in 8-move solitaire the central game is solvable
on all five boards.

To identify the holes on the board, we use a notation where the rows are
labeled, top to bottom, 1,2,3, \ldots,
and the columns are labeled left to right, a,b,c, \ldots,
as in Figure~\ref{fig1}b.
Note that in this notation the coordinate of the board center
is not always the same, the center is $d4$ for $Diamond(4)$
and the boards of Figure~\ref{fig1},
but is $c3$ for $Diamond(3)$ and $e5$ for $Diamond(5)$.

A solitaire jump is denoted by the starting and ending coordinates for the jump,
separated by a dash, i.e., d4-f4 for the rightward jump from the center in
Figure~\ref{fig2}b.
If the same peg makes one or more jumps, we will call this a \textbf{move}.
To denote this we add the intermediate coordinates of the jumps, for example
d4-f4-d6-d4 for the triple-jump move from the center in Figure~\ref{fig2}b.

A board position $B$ will be denoted by a capital letter, and
by $B'$ we mean the \textbf{complement} of this board position, where
every peg is replaced by a hole, and vice versa.
We will also refer to the board position with only
the hole $d4$ filled as $d4$; the central game on the
33-hole board is thus the problem of playing from $d4'$ to $d4$.
If $S$ is any subset of holes, and $B$ is a board position,
then $\#_{\mbox{\tiny\textit B}} S$ is the number
of holes in $S$ that are occupied by pegs.

A common type of peg solitaire problem begins with one peg
missing and finishes with one peg.
Beasley \cite{Beasley} refers to such problems as
\textbf{single vacancy to single survivor} problems,
or \textbf{SVSS} problems.
A \textbf{complement problem} is the special problem where the starting vacancy
and finishing hole are the same, because the starting and ending board
positions are complements of one another.
Complement problems are particularly attractive because of this symmetry,
and the central game is perhaps the most attractive of all because the
starting and ending board positions are square symmetric.

The solution to any SVSS problem on the standard 33-hole board has
exactly 31 jumps, because we begin with 32 pegs and finish with one,
and one peg is lost per jump.
However we can also consider the number of moves in a solution, which
has a natural interpretation as the number of pegs
that must be touched during the solution
(not counting those removed from the board).
A solution's \textbf{length} will always be measured in moves.
A classical problem in the history of 4-move solitaire is to
determine the shortest solution to all SVSS problems
on the standard 33-hole board.
In 1912, Ernest Bergholt \cite{Bergholt} found a solution
to the central game in 18~moves.
In 1964, John Beasley proved that the central game
can't be solved in under 18~moves \cite{Beasley, WinningWays},
so Bergholt's solution is the shortest possible.
Subsequently, shortest solutions to all SVSS problems on the
33-hole board have been found \cite{Beasley, WinningWays}
using computational search.

\vskip 30pt
\section*{\normalsize 2. Single vacancy to single survivor (SVSS) problems}

\subsection*{\normalsize 2.1 Reversibility, categories and position classes}

8-move and 4-move solitaire are similar in one very important way,
in regard to the \textbf{reversibility} of the game.
In 4-move solitaire, if a sequence of jumps takes one from
board position $A$ to $B$, then these jumps executed in
the same direction, but in reverse order will take one
from $B'$ to $A'$.
This stems from the fact that the basic solitaire jump itself
takes the complement of three consecutive board locations.
As stated in Winning Ways \cite{WinningWays},
``Backwards solitaire is just forward solitaire with the notions
``empty" and ``full" interchanged."
This property is not lost in 6-move or 8-move solitaire, and
these games also have the reversibility property.
One consequence of this is that any solution to a complement problem,
when played in reverse, is a different solution to the same
complement problem.

\begin{figure}[htbp]
\centering
\epsfig{file=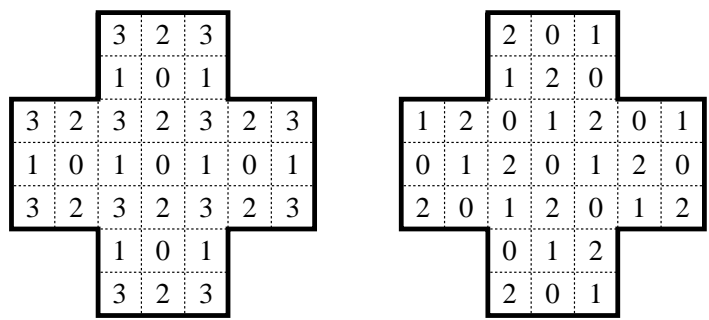}
\caption{(a) Four categories of pegs (b) The diagonal labeling for position classes}
\label{fig3}
\end{figure}

We can separate the pegs into four \textbf{categories} (0-3) by their
ability to jump into the center, jump horizontally over
the center, jump vertically over the center, or jump
diagonally over the center (Figure~\ref{fig3}a).
Any peg is in exactly one category, and remains in
this category for the entire game.
A peg cannot remove another peg of the same category,
and in 4-move solitaire a peg can remove pegs only in
two of the three other categories.
In 6-move and 8-move solitaire, a peg can remove any peg in
the other three categories (ignoring any limitation due
to the edge of the board).
This suggests that in general, the same SVSS problem can be
solved in fewer moves in 6- and 8-move solitaire.

Allowing diagonal jumps adds more potential jumps,
but also changes the game in a more fundamental way.
In 4-move solitaire, suppose we label the diagonals alternately
$0,1,2,0,1,...$ as in Figure~\ref{fig3}b,
and let $n_i$ be the number of pegs in the holes labeled $i$.
The key observation is that every (non-diagonal)
jump involves exactly one hole of each of the three labels.
After a jump is executed,
two of the three $n_i$ decrease by 1,
while the other increases by 1.
Therefore, if we add any pair among $\{n_0, n_1, n_2\}$,
the parity (even or odd) of this sum cannot change as the game is played.
There are six invariant parities---those
associated with the sums $n_0+n_1$, $n_1+n_2$, and $n_0+n_2$,
and the analogous sums along diagonals slanting up and right.
The set of all board positions can be partitioned into 16
\textbf{position classes} \cite{Beasley, Bell},
where all boards in a position class share identical values of the six parities.
During a game of 4-move solitaire,
the board position can never leave the starting position class.
The position classes can also be derived using
algebraic rules \cite{WinningWays}.

These position classes constrain the possible finishing
holes in a game of 4-move solitaire.
Using Cartesian coordinates,
all board positions with a single peg at $(x+3i,y+3j)$
are in the same position class for any integers $i$ and $j$.
This restriction is known as the ``Rule of Three" \cite{WinningWays}.
In the next section we'll see that in 8-move solitaire
there is no such restriction on finishing holes.

\subsection*{\normalsize 2.2 SVSS solvability in 8-move solitaire}

{\bf Theorem 1} Under 8-move solitaire, for all five boards except
$Diamond(3)$, any SVSS problem is solvable.
That is, beginning from a full board with one peg missing,
it is possible to finish with one peg at any board location.

\begin{figure}[htbp]
\centering
\epsfig{file=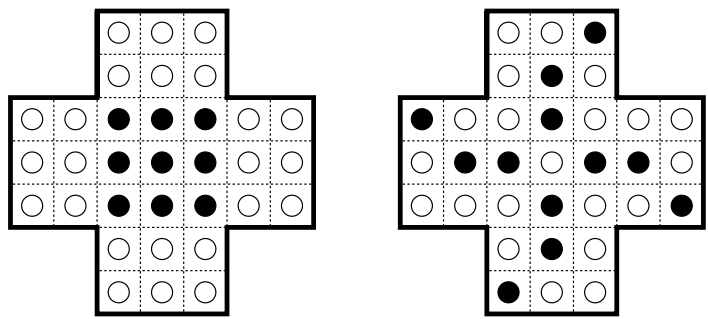}
\caption{The board position (a) C9, and (b) P12}
\label{fig5}
\end{figure}

{\it Proof:} Consider the board position of
Figure~\ref{fig5}a with only the central 9 board locations occupied.
We'll refer to this board position as ``$C9$".
In 4-move solitaire, $C9$ is in the position class of the empty board,
and it is therefore impossible to reach any single peg finish.
In 8-move solitaire, however,
it is possible to play from $C9$ to finish
with one peg anywhere on the board.
This is the case not only for the 33-hole board,
but also on the 37-hole board and $Diamond(5)$.
Because the board position $C9$ is square symmetric,
we need only list a few solutions to prove this---these solutions
can be found in Appendix A.

The proof will be complete if it is possible
to solve the $C9$~complement (play from $C9'$ to $C9$),
for we can then construct a solution to any
SVSS problem from $a'$ to $b$ as follows: (1) we take a solution
from C9 to $a$, and reverse it to obtain a solution from $a'$ to $C9'$,
(2) we then play from $C9'$ to $C9$ and finally
(3) from $C9$ to $b$.

The $C9$~complement problem is an interesting puzzle in itself,
and could be called the ``big central game."
The $C9$~complement is only solvable for the two largest boards,
the 37-hole board and $Diamond(5)$.
A solution to the problem on the 37-hole board is shown in Figure~\ref{fig4},
while a solution on $Diamond(5)$ is in Appendix A.
\begin{figure}[htbp]
\centering
\epsfig{file=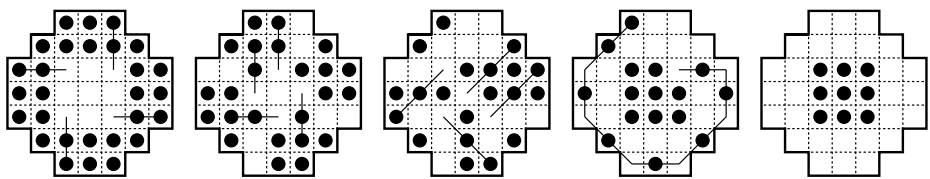}
\caption{A 13-move solution to the $C9$~complement on the 37-hole board.}
\label{fig4}
\end{figure}

It remains to prove the theorem in the case of the 33-hole board
and $Diamond(4)$.
In the first case we just need to work a bit harder.
Although the $C9$~complement is unsolvable (as will be proved shortly),
we can find another board position $P12$ (Figure~\ref{fig5}b)
from which any finishing location can also be reached (see Appendix A).
Moreover it is easy to play from $P12'$ to $C9$,
even without diagonal jumps.
Therefore our solution from $a'$ to $b$ goes from
$a'$ to $P12'$ to $C9$ to $b$.

For $Diamond(4)$, Theorem 1 does not fall to such elegant arguments,
and we have only demonstrated the result using an exhaustive search.
Fortunately, an exhaustive search goes fairly quickly since the
board has only 25 holes.

\begin{figure}[htbp]
\centering
\epsfig{file=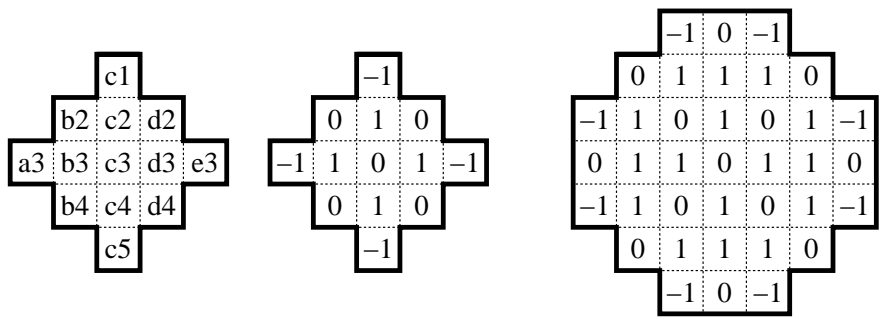}
\caption{(a) $Diamond(3)$ coordinates, (b) resource count and
(c) a resource count on the 33 or 37-hole boards.}
\label{fig6}
\end{figure}

For $Diamond(3)$\footnote{The $Diamond(3)$ board can be purchased
under the name ``Hoppers," marketed by ThinkFun\texttrademark.},
we can easily prove that not all SVSS problems
are solvable using the \textbf{resource count} or
\textbf{pagoda function} in Figure~\ref{fig6}b.
A resource count is a weighting function on the board,
carefully devised so that no jump increases its total.
To calculate the value of the resource count for a particular
board position we simply sum the numbers where a peg is present.
If we begin a SVSS problem with a vacancy at one of the holes
marked ``1" in Figure~\ref{fig6}b, then the resource count
starts at -1,
and it is therefore impossible to finish with one peg
at any board location with a weight greater than -1.

These restrictions are only necessary conditions, however, and
there are SVSS problems on $Diamond(3)$ feasible by the
resource count that are not solvable.
An exhaustive search demonstrates, for example, that if we start with
a vacancy at $c2$, we can finish only at $c1$ or $c5$,
not at $a3$ or $e3$.

The reader is cautioned that the majority of
resource counts for 4-move solitaire are not valid in
8-move solitaire.
Careful checking of all jumps, regular and diagonal,
is required to ensure that a particular weighting is
a valid resource count.

The ``cross count" of Figure~\ref{fig6}c is a valid
resource count for 8-move solitaire on the 33 or 37-hole boards.
For the $C9$~complement, the value of the resource begins at 4
and ends at 4, so a solution cannot include any jump that
reduces this resource count.
One such jump is a vertical or horizontal jump into the center;
this jump loses 2 and is the only way on the 33-hole to fill the
central hole, hence the $C9$~complement is unsolvable on this board.
On the 37-hole board, we see that the central hole must be 
filled by a diagonal jump.

In fact a similar argument using the resource count of
Figure~\ref{fig6}c proves that on the 33-hole board
there is no $90^\circ$ rotationally symmetric board position $A$ that
satisfies both:
(1) any finishing hole can be reached from $A$ and
(2) It is possible to play from $A'$ to $A$.

\vskip 30pt
\section*{\normalsize 3. Short solutions}

John Beasley \cite[p. 232]{Beasley} gives a 16-move solution
to the central game on the 33-hole board under 8-move solitaire.
His comment ``but I suspect that this can be improved"
was the impetus for this paper.
What is the shortest solution to the central game in 8-move solitaire?
In this section we will answer this question for
all five board types.
In fact we will also answer a more general question:
what is the shortest solution to any SVSS problem on this board?

\subsection*{\normalsize 3.1 Computational search techniques}

On a board with $n$ holes, the total number of possible
board positions is $2^n$; this is an upper bound on the
size of the state space for the game\footnote{In 4-move solitaire, this upper bound can be reduced to
the size of the starting position class, which is $16$ times smaller, or $2^{n-4}$.}.
For 8-move solitaire on $Diamond(4)$,
the upper bound is $2^{25}\approx 3.4\times 10^7$.
From $d4'$, the number of board positions that can be reached
can be calculated as $2.7\times 10^7$,
about 80\% of this upper bound.
For larger boards in 8-move solitaire,
it is a reasonable estimate that 70-90\% of the upper bound
can be reached from any single vacancy start,
which for the 37-hole board gives an estimated state space
of $1.1\times 10^{11}$ nodes.
Of course, the problem starting with $d4$ missing has 8-fold symmetry,
and this can be used to reduce all these estimates by
nearly a factor of 8.

The computational technique used to find short solutions in this paper is
called a search by levels \cite[p. 249]{Beasley}.
It is a breadth-first search by moves.
Each level set $L_i$ is a set of board positions,
$L_0$ is the set containing only the initial board position.
To calculate $L_i$ from $L_{i-1}$, we take all board positions in $L_{i-1}$,
and execute all possible single moves.
After removing duplicates and board positions seen at previous levels, we
obtain $L_i$, the set of board positions reachable
from the starting position in exactly $i$~moves.
If the target board position appears in some $L_m$, then
this problem can be solved in a minimum of $m$ moves.
Moreover, by tracing the target board
position back through the level sets,
we can find \textbf{all} solutions of length $m$.

The reason a breadth-first rather than a depth-first search is
used to find short solutions in peg solitaire is
the amount of repetition encountered while searching.
For example, if one can play from board position $A$ to $B$ in 8~moves,
it is common to have millions of possible
move sequences that can take one from $A$ to $B$.
To effectively apply a depth-first search, it is necessary
to store board positions encountered previously, and this
can fill memory for the largest boards.
This memory limitation is also present in a breadth-first search,
but we have the option to store only
board positions at the current level,
rather than all previously seen.
A breadth-first search in peg solitaire can also be easily split
into smaller pieces to be worked on separately,
and we can also eliminate searches over boards that
are equivalent by symmetry.

One search technique that has been successfully applied
to find short solutions in peg solitaire is
a bidirectional search \cite{BGS}.
Here a search by levels is run forward from the initial board state,
then backwards from the final board state,
with a solution identified in the intersection of the two searches.
While useful, we'll see that this search technique has
a fundamental limitation when applied to larger peg
solitaire boards, particularly in 8-move solitaire.

A second technique, the A* search \cite{Astar},
uses an estimate $h(B)$ of the number of moves from
the current board position $B$ to the desired goal
(usually a one peg finish).
A* search techniques have been successfully applied to other puzzles,
such as the Fifteen Puzzle \cite{Korf85,Korf96}.
To apply an A* search in our search by levels, at each level $i$,
we accept only board positions satisfying the constraint
\begin{equation}
i+h(B)\le m,
\label{astar}
\end{equation}
where $m$ is the length of the longest solution
we will accept.
This is not a traditional A* search, in which the node
selected for expansion is the one with the
smallest value of $i+h(B)$.

To find the shortest solution to a problem starting
from board position $S$,
we set $m=h(S)$ and run a search by levels with
constraint (\ref{astar}).
If this finishes with no solution found,
we set $m=h(S)+1$ and repeat the search from the start,
and continue this process, increasing $m$ by one each time,
until a solution is found.
This technique is very similar to
``iterative deepening A*" \cite{Korf85},
the main difference being that it is based on
a breadth-first rather than a depth-first search.
Although we haven't tried a depth-first search,
it seems likely that this technique will find
a solution more quickly than a breadth-first
search by levels.
One advantage of our technique is that we can
find all solutions, rather than just the first one.

If $h(B)$ is \textit{admissible},
meaning that it never overestimates the actual number
of moves remaining, then an A* search is guaranteed to
find the shortest solution, if one exists \cite{Astar}.
For this reason we will always select $h(B)$
that are admissible.
One requirement of admissibility is that if $T$ is a
target board position, $h(T)=0$.

\subsection*{\normalsize 3.2 Corner constrained boards}  

We call a board location a \textbf{corner} if there is no
jump that can capture a peg at this board location.
Note that whether a board location is a corner depends
upon the particular type of solitaire---for example in
8-move solitaire $Diamond(n)$ has 4 corners,
but in 4-move solitaire every hole at the edge of the board is a corner.
A \textbf{corner peg} is defined as any peg that is in
the same category as some corner peg, or equivalently
the peg can jump into some corner.
Finally, a \textbf{corner move} is any move that begins from a corner.

Let us now consider $Diamond(3)$.
This board has four corners ($c1$, $a3$, $e3$, and $c5$),
and the corner pegs are those in the corners, plus the center $c3$.
A peg which begins in a corner cannot be removed in its
original location, but must be first moved to the center.
Notice that no move can remove more than one corner peg,
the central peg.
Thus the process of removing each peg that begins in a corner
involves at least two moves: moving it to the center and then
jumping over it.

If we consider the central game on $Diamond(3)$, we must
move three of the corner pegs to the center and then
jump them, one at a time, and on the final move bring
the fourth corner peg to the center.
Therefore, no solution to the central game can have fewer than 7~moves.
Since a solution can be found in 7~moves (Appendix A),
we have proved that it is of minimum length.

Similar arguments work on any board that has the following 
(qualitative) properties:\begin{packed_enumerate}
\item All corner pegs are in the same category
(which implies that no corner move can remove a corner peg).
\item There is a limit to the number of corner pegs that can be removed by one move.
Strictly speaking, of course, there is always some limit
to the number of pegs removed by one move, but this property
is reserved for boards for which this limit is small enough
to constrain the length of solutions.
\end{packed_enumerate}
A board with the above two properties is called
\textbf{corner constrained}, because the length of the shortest solution
is limited primarily by the removal of corner pegs.
Of the five boards we consider, the corner constrained boards are
$Diamond(3)$, the 33-hole board and the 37-hole board.
The other two boards are called \textbf{edge constrained} and finding
short solutions requires a different technique.

Let us consider now the 37-hole board, because much can be proved
about this board, without resorting to a computational search.

{\bf Theorem 2} Under 8-move solitaire on the 37-hole board:
\begin{packed_enumerate}
\item The central game cannot be solved in less than 13~moves.
\item No SVSS problem can be solved in less than 12~moves.
\item The $C9$~complement cannot be solved in less than 13~moves,
and all the moves of a 13-move solution must start outside
and end inside the central 9 holes.
\end{packed_enumerate}

{\it Proof:} This board has 8 corners and the central game begins
with 12 corner pegs and finishes with zero.
In addition, no move can remove more than four corner pegs
(those at $c3$, $e3$, $c5$ and $e5$).
Therefore we will need 8 corner moves,
plus at least 3 more moves to remove
the 12 corner pegs, for a total of 11~moves.
However, the first move can only be d2-d4 or b2-d4
(or symmetric equivalents), which remove at most one
corner peg, so we need at least 12~moves.

But it is not hard to see that even 12~moves is impossible.
For this requires that after the first move, every move either
(1) is a corner move that does not finish at a corner, or
(2) removes at least 3 corner pegs.
If the first move is d2-d4, then there is no second move of
either kind.
If the first move is b2-d4, then we can play c1-c3 or e1-c3,
but then we can no longer make a move of either kind.

For the second part of Theorem 2 we use a similar argument.
For a general SVSS problem, we begin with at least 11 corner pegs
and finish with one or zero, so we must remove at least 10.
In addition there must be at least one move beginning from every corner.
This is clear if the corners are all filled at the start, and even
if the starting vacancy is at a corner, it is filled by the first move,
and therefore there must still be a move out of it.
This gives a total of 11~moves, but 11~moves cannot be attained.
If we begin anywhere but a corner, then every move must either
be a corner move that does not finish at a corner,
or remove two corner pegs,
but there is no way that the first two moves can be like this.
If we begin at a corner, then the first move is arbitrary but
there is no way the second and third moves can be like this.

For the third part of Theorem 2, the $C9$ complement,
we must again have eight corner moves.
These eight moves can at best leave pegs
at the four holes $c3$, $e3$, $c5$ and $e5$, which leaves five holes in the center
that must be filled by at least five more moves, for a total of 13 moves.
Note that the solution in Figure~\ref{fig4} has minimum length.
For any 13-move solution, nine of the moves must supply a peg to the central
region, and the other four moves remove a corner peg and cannot reduce the number
of pegs in the central region (otherwise we will need more than 13~moves).
Therefore all the moves must start outside and end inside this region.

The ideas of the above proof can also be incorporated into
numerical schemes to find a 13-move
solution to the central game,
a 12-move solution to the $c1$ or $c3$ complements,
and a 13-move solution to the $C9$~complement.
For any board position $B$ we come up with a lower bound
$h_1(B)$ on the number of moves to any single peg finish.
Let $c(B)$ be the number of corners occupied by pegs, and
$p(B)$ be the number of corner pegs.
In our notation, $c(B)=\#_{\mbox{\tiny B}}\{c1,e1,a3,g3,a5,g5,c7,e7\}$ and
$p(B)=c(B)+\#_{\mbox{\tiny B}}\{c3,e3,c5,e5\}$.
We must have a move starting from each occupied corner, and separately
remove the corner pegs at most 4 per move.
So the number of moves from $B$ to a single peg finish is at least
\begin{equation}
h_1(B)=c(B)+\lceil(p(B)-f)/4\rceil,  \label{hdef1}
\end{equation}
where the ceiling operator specifies rounding up to the nearest integer,
and $f$ is set to either 0 or 1 depending on which finishing
locations we allow.
Setting $f=0$ will allow any single peg finish except for a corner peg.
If we want to find solutions that finish with one peg anywhere,
we use $f=1$, and if $B$ has exactly one peg, we override (\ref{hdef1})
and set $h_1(B)=0$.
Note that the heuristic (\ref{hdef1}) is valid only for the
33- and 37-hole boards.

\begin{figure}[htbp]
\centering
\epsfig{file=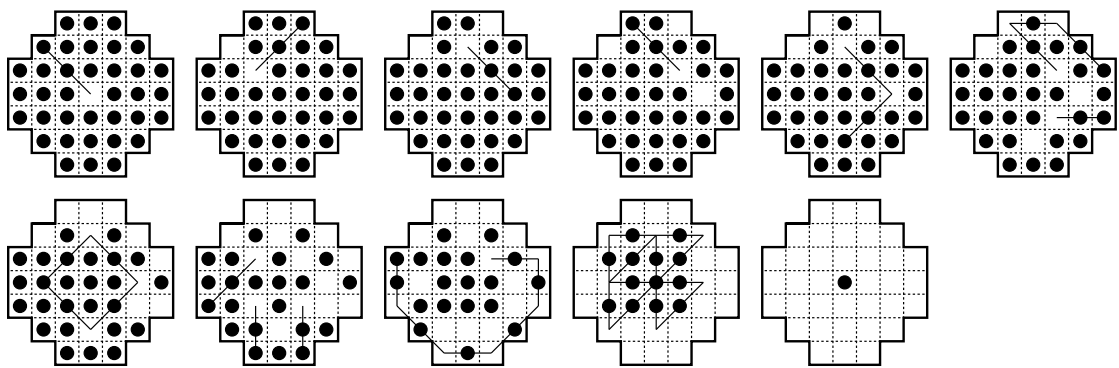}
\caption{A minimum length, 13-move solution to the central game on the 37-hole board.}
\label{fig7}
\end{figure}

Using the heuristic (\ref{hdef1}) in an A* search (\ref{astar})
we can find all 13-move solutions to the central game in a matter
of minutes (Figure~\ref{fig7}).
Our search algorithm also takes advantage of the 8-fold
symmetry of the problem,
and expands a total of $2.9\times 10^6$ nodes,
a reduction by a factor of $4700$ over an unconstrained search
(using the previous estimate of $1.1\times 10^{11}/8$
for the total number of nodes).

Note that after the third move in Figure~\ref{fig7},
constraint (\ref{astar}) is satisfied with an equality,
and this is the case for nearly all nodes in the search.
Normally in a search by levels, one needs to eliminate
boards at the current level that have appeared in previous levels.
However when most nodes satisfy the A* constraint (\ref{astar})
with an equality, this check is no longer necessary.
We can save memory because we need only to keep
boards at the current level.

One subtlety of this search is that it takes place
first by individual moves and second
by extending the current move.
It is important that the constraint (\ref{astar}) be applied only
after the search over the current move is finished.
For example, if we examine any board position of Figure~\ref{fig7}
in the middle of a multi-jump move, the constraint (\ref{astar})
is violated.

The heuristic (\ref{hdef1}) is based on finding solutions with single
peg finishes, so a different $h(B)$ is needed for the
$C9$~complement problem.
We still need one move per corner peg, and these moves can only fill
four of the nine vacancies in the central region.
So we can add at least one move for each of the remaining 5 vacancies.
If we let $v(B)=5-\#_{\mbox{\tiny B}}\{d3,c4,d4,e4,d5\}$,
then we use the heuristic
\begin{equation}
h_2(B)=c(B)+v(B),  \label{hdef2}
\end{equation}
which finds all 13-move solutions in less than 1 minute after
expanding $2.0\times 10^5$ nodes.
Our algorithm finds that all 13-move solutions to the $C9$
complement contain exactly 7 diagonal jumps.

\begin{figure}[htbp]
\centering
\epsfig{file=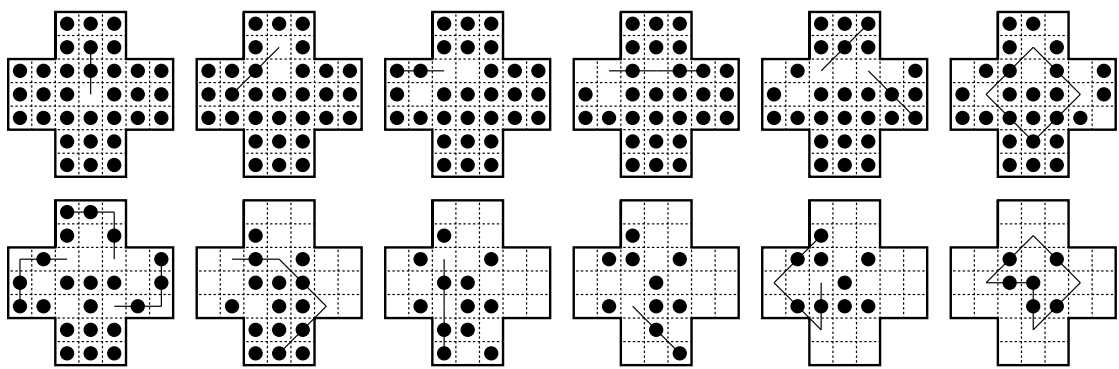}
\caption{A minimum length, 15-move solution to the central game on the 33-hole board.}
\label{fig8}
\end{figure}
The arguments of Theorem~2 can be applied to the standard 33-hole board to
prove that there is no solution to the central game with less than 14~moves.
However, an A* search (\ref{astar}) finds that the
shortest solution has 15~moves (Figure~\ref{fig8}), after
expanding $1.3\times 10^7$ nodes, a speed increase
of a factor of $70$ over an unconstrained search.
The solution in Figure~\ref{fig8} is identical to the 16-move
solution given by Beasley \cite[p. 241]{Beasley} up until the 9th move.
Using an A* search, we can also find 13-move solutions
to the $c1$ and $c3$ complements (Appendix A),
and that no SVSS problem on this board is solvable in less than 13~moves.

\subsection*{\normalsize 3.3 Edge constrained boards}

The corner pegs on $Diamond(4)$ are not all in the same category,
and on $Diamond(5)$ a single move can remove 9 corner pegs,
so the techniques of the previous section do not apply.
However, using a different argument we can see that the
presence of the edge of the board (including the corners) significantly
constrains the solution length.
\begin{figure}[htbp]
\centering
\epsfig{file=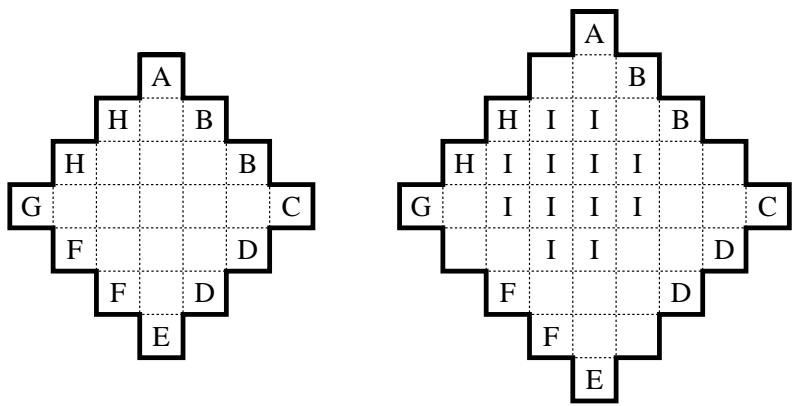}
\caption{$Diamond(4)$ and $Diamond(5)$ divided into ``Merson regions," A,B,C,\ldots.}
\label{fig9}
\end{figure}

Consider $Diamond(4)$, and divide the board into eight
``Merson regions\footnote{Named after Robin Merson who first used this concept in 1962 on the 6x6 square board \cite[p. 203]{Beasley}.}"
A-H as in Figure~\ref{fig9}a.
The shape of a region is chosen such that when it is entirely filled with pegs,
there is no way to remove a peg in the region without a move that originates
inside the region.

{\bf Theorem 3} Under 8-move solitaire on $Diamond(4)$:
\begin{packed_enumerate}
\item The central game cannot be solved in less than 10~moves.
\item No SVSS problem can be solved in less than 8~moves.
\end{packed_enumerate}

{\it Proof:} At the start of the central game, we have 8 regions full,
and neither the first nor the last move can start from these regions,
so we need at least 10 moves.
Solutions with 10 moves can be found (Appendix A), so they have minimum length.

By the same argument, if we begin anywhere in the interior
of the board, any SVSS problem requires at least 8 moves.
If we begin at a corner, this region is filled by the first move,
and we will still need at least 7 more moves.
The only remaining possibility for a 7 move SVSS solution
is to start along an edge, say $c2$.
Every jump must come from the 7 filled regions (A-G), and we cannot
jump into a corner or refill any edge region.
This forces the jumps: a4-c2, d1-b3, e2-c2, g4-e2, which
leaves us in a board position where there is no move meeting
these constraints.
So no SVSS problem is solvable in 7 moves.
Yet 8 move solutions exist (see Appendix A).

These ideas can be used in an A* search to
find short solutions on $Diamond(5)$.
From each board position $B$, we again obtain a lower bound on the
number of moves remaining to any single peg finish.
Referring to Figure~\ref{fig9}b for $Diamond(5)$, we define
\begin{equation}
h_3(B)=\mbox{the number of filled regions.}
\label{hdef3}
\end{equation}
This heuristic can be improved by ``sliding" the edge or interior regions
to any configuration where they do not overlap.
For the best bound, we initialize $h_3$ to zero,
then add +1 for every occupied corner, and for every edge
which has two consecutive holes filled (not including the corners).
We then add +1 if region ``I" is completely filled by pegs,
either in the configuration in Figure~\ref{fig9}b, or
by translating region ``I" one hole to the right, down, or both.
Again, if $B$ has only one peg then we override (\ref{hdef3})
and set $h_3(B)=0$.

Another algorithm improvement comes from the following idea:
if, for some board position $B$ at level $i$, we have $i+h_3(B)=m$,
then \textit{all remaining moves must begin from filled regions}.
For example, the interior region ``I" is unlikely
to be filled after the first few moves.
For any board position $B$ with ``I" unfilled and $i+h_3(B)=m$,
all remaining moves must start from the edge of the board (including corners).

\begin{table}[htbp]
\begin{center} 
\begin{tabular}{ | c | r | r | r | }
\hline
 & \multicolumn{3}{|c|}{Size of level set, $|L_i|$} \\		
Level \# ($i$) & unconstrained & $m=10$ & $m=11$  \\
\hline
\hline
0 & 1 & 1 & 1 \\
\hline
1 & 2 & 2 & 2 \\
\hline
2 & 12 & 10 & 12 \\
\hline
3 & 152 & 66 & 139 \\
\hline
4 & 2,347 & 216 & 1,381 \\
\hline
5 & 43,763 & 630 & 9,134 \\
\hline
6 & 890,355 & 2,002 & 54,798 \\
\hline
7 & 18,085,322 & 6,007 & 372,122 \\
\hline
8 & 325,165,209 & 17,497 & 2,739,963 \\
\hline
9 & not calculated & 2,637 & 20,776,877 \\
\hline
10 & not calculated & 0 & 15,467,734 \\
\hline
11 & not calculated & 0 & 8 \\
\hline
Total nodes expanded & est $2.2\times 10^{11}$ & $2.9\times 10^4$ & $3.9\times 10^7$ \\
\hline
Total run time & & 0.8 min & 10 hrs \\
\hline
\end{tabular}
\caption{Results from a search by levels for the central game on $Diamond(5)$.} 
\label{table1}
\end{center} 
\end{table}

Table~\ref{table1} shows a result of a search for the shortest solution
to the central game on $Diamond(5)$ with and without the
A* constraint $i+h_3(B)\le m$.
Here we have included the 8-fold symmetry of the problem
(there are only 2 different first moves).
The second column shows the size of the level sets when the search is
unconstrained.
Note that to store the level set $L_8$ requires over 1.8 gigabytes of disk space,
using 6 bytes per board position,
and this search cannot be run to completion due to
excessive time and disk space requirements.
In the third column, we use constraint (\ref{astar}) with $m=10$,
which shows that no solution of length 10 exists
(ending at the center or any other board location).
To find an 11-move solution, we can begin from $L_8$
in the $m=10$ column, but relax the constraint to $m=11$.
To find \textbf{all} 11-move solutions, however,
we must start over from scratch with $m=11$ as shown
in the final column.

In the final two columns of table~\ref{table1}, note that
$|L_i|$ decreases at $i=m-1$.
The reason for this is an improvement in the $h_3(B)$ constraint
when $i=m-1$.
For this level only,
let $r$ be the number of filled regions as in (\ref{hdef3}), and
$$
s=\left\{ \begin{array}{rl}
0 & \mbox{if there is only one peg remaining,} \\
1 & \mbox{if there is exactly one peg in at least one category,} \\
2 & \mbox{otherwise.}
\end{array} \right.
$$
Then we use $h_3(B)=\max(r,s)$.
This is valid because before the final move,
the finishing peg must be the only peg in its category.

In the final column, $|L_{11}|=8$ indicates that there are
8 possible finishing locations.
Because of the symmetry of the problem,
these 8 locations cover most of the board,
and this calculation shows that,
beginning from $e5'$,
there exists an 11-move solution finishing anywhere except
the 4 corners.
See Appendix~A for an 11-move solution to
the central game.

Can we find an 11-move solution to the central game
more quickly using a bidirectional search?
Interestingly, a bidirectional search performs poorly on this problem.
To compute one level backwards from the final board position,
the algorithm must calculate all board positions that
can finish in one move to $e5$.
There are over half a million such board positions,
even counting symmetrical board positions as the same.
The fact that the search tree grows so rapidly is
not the only problem,
because most of these final moves are long
sequences of nested loops up to 24 jumps long
(for example see the final move in Figure~\ref{fig11}).
A move with $p$ loops can be executed at least $2^p$ ways
(each loop can be traversed in either direction),
so the number of ways to traverse each move grows
exponentially with $p$.
Our program takes almost 9 hours to go back one level,
and in 24 hours has completed less than 1\% of the
calculation back two levels.

\begin{figure}[htbp]
\centering
\epsfig{file=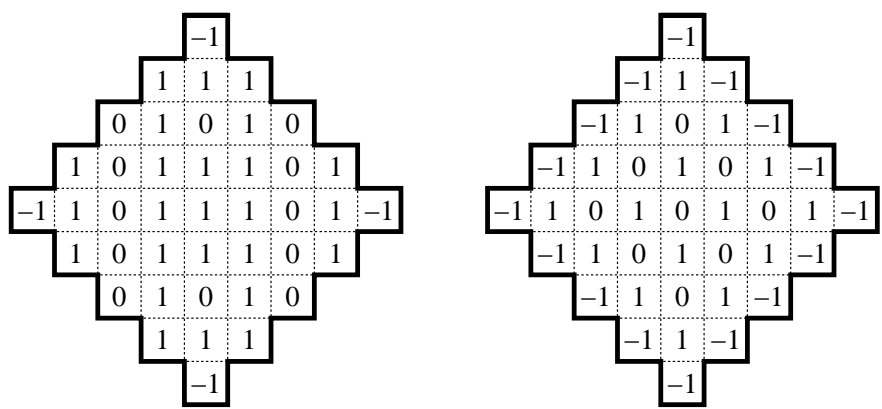}
\caption{Resource counts on $Diamond(5)$ useful for finding short solutions to (a) the $C9$~complement in 8-move solitaire and (b) SVSS problems in 4-move solitaire.}
\label{fig10}
\end{figure}

Finding the shortest solution to the $C9$ complement on $Diamond(5)$ is
difficult.
The resource count of Figure~\ref{fig10}a is useful for this computation,
because for the $C9$~complement it starts at 12 and ends at 9, so
we can afford to lose only 3.
We can apply this in our search by adding the constraint that
this resource count must be at least nine.
We can also use the A* search heuristic (\ref{hdef2}) with
$c(B)=\#_{\mbox{\tiny B}}\{e1,a5,i5,e9\}$ and
$v(B)=8-\#_{\mbox{\tiny B}}\{d4,e4,f4,d5,f5,d6,e6,f6\}$.
However, even with both these constraints,
an A* search to level 16 or 17 is time consuming.
A bidirectional search is the fastest technique for this problem,
since the last move cannot be a long sequence of loops.
We ran forward 9 levels and backward 7, and found the
intersection between these level sets was empty,
so the $C9$~complement cannot be solved in 16~moves.
A 17-move solution can be found in Appendix A.

\subsection*{\normalsize 3.4 Back to 4-move solitaire}

The A* heuristic (\ref{astar}) has proved so useful in 8-move solitaire,
why not use it to find short solutions in 4-move solitaire?
On the standard 33-hole board, we can indeed do so,
and the heuristic (\ref{hdef1}) finds Bergholt's 18-move solution
to the central game in under a minute, after expanding $4.6\times 10^5$
nodes\footnote{The most efficient way to do this is to notice that any 18-move
solution from $d4'$ to $d4$ generates a 17-move solution from $d4'$ to $d1$
(and vice versa).
Thus, in (\ref{astar}) we set $m=17$, and in (\ref{hdef1}) $f=0$.}.
This is a reduction by about a factor of $50$
over an unconstrained search, for which the total number of
reachable board positions is $2.3\times 10^7$.

The 37-hole and $Diamond(n)$ boards have additional corners
in 4-move solitaire,
and although the heuristics (\ref{hdef1}) or (\ref{hdef3})
are still valid they hardly speed up the search at all.
On the 37-hole board under 4-move solitaire, there are only
10 solvable SVSS problems, and an exhaustive
search has found that all require a minimum of
20 or 21~moves \cite{GPJ28}.
The resource count in Figure~\ref{fig10}b is valid on this board,
and can speed up the search by a factor of two or so.
On the 37-hole board, it actually takes less
time to find short solutions to SVSS problems under
8-move solitaire than 4-move solitaire.

On $Diamond(5)$ under 4-move solitaire, there are only four SVSS problems
solvable, all in a minimum of 26~moves \cite{BGS}.
The resource count in Figure~\ref{fig10}b is very useful in speeding
up this search, but again finding minimum length
solutions to SVSS problems is comparable or easier under 8-move solitaire.

\vskip 30pt
\section*{\normalsize 4. Conclusions}

We have found that 8-move solitaire is similar to 4-move solitaire,
except for the major difference that, in general,
any SVSS problem is solvable.
Table~\ref{table2} summarizes our results on short solutions to
SVSS problems for the five boards.
Table~\ref{table2} shows that as boards become larger, rather
counter-intuitively the length of the shortest solution may decrease!
In going from the 33-hole board to the 37-hole board to $Diamond(5)$,
4 holes are added each time, yet the length of the shortest solution
to the central game \textit{decreases} from 15 to 11.

\begin{table}[htbp]
\begin{center} 
\begin{tabular}{ | l | c | l | c | c | }
\hline
 & Holes / & & \multicolumn{2}{| c |}{shortest solution in moves to} \\		
Board & Corners & Type & the central game & any SVSS problem  \\
\hline
\hline
$Diamond(3)$ & 13 / 4 & corner constrained & 7$\dagger$ & 7$\dagger$ \\
\hline
$Diamond(4)$ & 25 / 4 & edge constrained & 10$\dagger$ & 8$\dagger$ \\
\hline
Standard 33-hole & 33 / 8 & corner constrained & 15 & 13 \\
\hline
37-hole & 37 / 8 & corner constrained & 13$\dagger$ & 12$\dagger$ \\
\hline
$Diamond(5)$ & 41 / 4 & edge constrained & 11 & 11 \\
\hline
\end{tabular}
\caption{Shortest solutions in 8-move solitaire ($\dagger$ - proved analytically).} 
\label{table2}
\end{center} 
\end{table}

We have discussed a number of techniques that can speed up
the search for the shortest solution,
sometimes by several orders of magnitude:
\begin{packed_enumerate}
\item A modified A* search (\ref{astar}) using an estimate
of the number of remaining moves, given by (\ref{hdef1}),
(\ref{hdef2}) or (\ref{hdef3}).
\item A bidirectional search.
\item A constraint from a resource count or pagoda function.
\end{packed_enumerate}

There is no single technique which works well for all peg solitaire
boards and all types of peg solitaire.
The A* search heuristic seems best suited in 8-move solitaire,
where it can speed searches up by a factor of 5000.
A bidirectional search works best in 4-move solitaire,
but fails on large boards where the final move
can be a long sequence of loops.
A constraint from a resource count often helps,
but its usefulness varies with the board and problem.
A resource count constraint can be applied at any time during a solution,
and is often effective in combination with one of the other
two techniques.

For many peg solitaire problems,
the constraint on filled Merson regions (\ref{hdef3})
is the most generally useful.
We have applied this technique successfully on triangular boards in
6-move solitaire to find the shortest solution
to SVSS problems on boards with up to 55 holes.

The game of 8-move solitaire probably will never be popular
because very complicated moves are possible.
Finding short solutions by hand seems virtually impossible.
One puzzle that is interesting to solve by hand is
the $C9$~complement or ``big central game" on the
37-hole board or $Diamond(5)$ (without minimizing moves).
This puzzle looks easy but is harder than it appears.

Given enough computational power, it is possible to find minimum length
solutions on all of these boards by exhaustive search.
However, in this paper we have tried to avoid a brute-force approach,
always looking for clever search heuristics to find solutions quickly.
Using an A* search, it can be faster to find minimum length solutions
on the same board under 8-move solitaire than under 4-move solitaire.
This is somewhat surprising since the number of board positions reachable
is about 16 times larger.

\vskip 15pt
\section*{\normalsize Acknowledgment}
The author would like to thank John Beasley for helpful discussions,
and the anonymous referee for suggestions regarding heuristic search techniques.

\vskip 30pt
\section*{\normalsize Appendix A. Solutions}

\subsection*{\normalsize A.1 Solutions from $C9$ to any finishing hole}

Because the 33-hole board and the board position $C9$ are square symmetric, it
suffices to list solutions to a single peg finish at
$d1$, $d2$, $d3$, $d4$, $c1$, $c2$ and $c3$.
On the 37-hole board, we also need to show that $b2$ can be reached, and
on $Diamond(5)$, the corner $e1$.

To d1: d4-b4-d6-f4-d2, c3-e3, d5-f3-d3-d1;
to d2: d4-f4-d6, c5-e5, d6-f4-d2-b4-d4-d2;
to d3: d4-d2-f4-d6-b4-d4, e4-c4, d5-b3-d3;
to d4: d4-b4-d6-f4-d2, c3-e3, d5-f3-d3, d2-d4;
to c1: d4-f4-d6, c4-c6-e6-c4-c2, e3-c3-c1;
to c2: d4-b4-d6, e4-e6-c6-e4-e2-c4-c2;
to c3: d4-b4-d2, e4-e2-c2-e4-e6-c4, c5-c3.

To b2 (37-hole board): e4-c2, d5-b3, d4-b2-d2-f4-d6-b4-b2.

To e1 (41 hole board): d5-d3, f4-d4-d2, f5-f7-d5, d6-d4, e5-c3-e1
(note the coordinate change).

\subsection*{\normalsize A.2 Solutions from $P12$ to any finishing hole}

In this case we consider only the 33-hole board.
The board position $P12$ is rotationally symmetric.
First we show how to reach the holes $d1$, $d2$, $d3$, $d4$ or any symmetric
equivalent.

To d1: a3-c5, g5-e3, d3-f3, c7-e5-e3, c4-c6-e4-e2, e1-e3, f3-d3-d1;
to d2: c7-e5, a3-c5, e1-c3-e3, d5-f5, g5-e5, f4-d6-b4-d4-f4-d2;
to d3: e1-c3, a3-c5-e5, d3-b3-d5-f5, c7-e5, g5-e3, f5-d5-f3-d3;
to d4: a3-c5, g5-e3, d3-f3, c7-e5-e3, c4-c6-e4-e2, e1-e3, f3-d3, d2-d4.

To get to any other finishing hole, we play e1-c3, a3-c5, c7-e5, g5-e3; then
to c1: e4-e6, c4-c6-e4, e3-e5, e6-e4-c2, c3-c1;
to c2: c4-e6, d3-f3-d5, e6-e4-c6-c4-c2;
to c3: c4-e2, d5-f5-d3, e2-e4-c2-c4, c5-c3.

\subsection*{\normalsize A.3 Minimum length solutions}

These solutions were found using a C++ program running on a PC of modest
speed (1GHz clock speed with 512MB of RAM).
Diagrams of most of these solutions can be found at\newline
{\tt http://www.geocities.com/gibell.geo/pegsolitaire/diagonal}

\noindent\textbf{$Diamond(3)$ solution}\newline
Central game in 7~moves\footnote{This is the same solution given by Beasley \cite[p. 241]{Beasley}.}:
c1-c3, c4-c2, a3-c3, d3-b3, c5-a3-c1-c3, d2-b4, e3-c5-a3-c3.

\noindent\textbf{$Diamond(4)$ solutions}\newline
Central game in 10~moves:
d2-d4, b5-d3, e2-c4, g4-e2, d1-f3, e6-g4-e2, b3-d1-f3-f5-d3-b5-b3-d3, d7-b5, a4-c6-e6-c4, d4-b4-d6-f4-d2-d4.
$f3$ complement in 8~moves:
d1-f3, g4-e2, d3-d1-f3, e6-g4-e2, b5-d3-f5, d7-b5, a4-c6-e6-e4-g4-e6-c4-e4, b3-d3-f3-d1-b3-b5-d5-f3.

\noindent\textbf{33-hole board solutions}\newline
Central game in 15~moves (Figure~\ref{fig8}):
d2-d4, b4-d2, a3-c3, f3-d3-b3, e1-c3, g5-e3, d6-b4-d2-f4-d6, c1-e1-e3, g3-g5-e5, a5-a3-c3, d7-f5-d3-b3, c7-c5-c3, e7-c5, c2-a4-c6-c4, d4-d6-f4-d2-b4-d4.
$c3$ complement in 13~moves:
a3-c3, d3-b3, f3-d3, e1-c3-e3, d1-f3-d3, g5-e3, a5-c3, d6-b4-d2-f4-d6, c1-c3-e3, c7-c5-c3, e7-c5, d4-f4-d2-b4-d6, g3-g5-e5-e7-c7-e5-c5-a5-a3-c3.

\noindent\textbf{37-hole board solutions}\newline
Central game in 13~moves (Figure~\ref{fig7}):
b2-d4, e1-c3, f4-d2, c1-e3, d6-f4-d2, g5-e5, g3-e1-c1-e3, b4-d6-f4-d2-b4, c7-c5, e7-e5, a5-c3, a3-a5-c7-e7-g5-g3-e3, d4-f2-d2-b2-b4-b6-d4-f4-d6-d4-d2-b4-d4.
$c3$ complement in 12~moves:
e1-c3, f4-d2, c1-e3, d6-f4-d2, a3-c1-e3, g5-e5, b4-d2-f4-d6-b4, c7-c5, e7-e5, g3-e3, b4-d6-f4-d2, a5-c7-e7-g5-g3-e1-e3-e5-c5-c3-e3-c5-a5-a3-c3-c1-e1-c3.
$C9$~complement in 13~moves (Figure~\ref{fig4}):
c7-c5, g5-e5, e1-e3, a3-c3, e6-e4, b5-d5, c2-c4, d1-d3, e7-c5, f2-d4, a5-c3, g3-e5, c1-a3-a5-c7-e7-g5-g3-e3.

\begin{figure}[htbp]
\centering
\epsfig{file=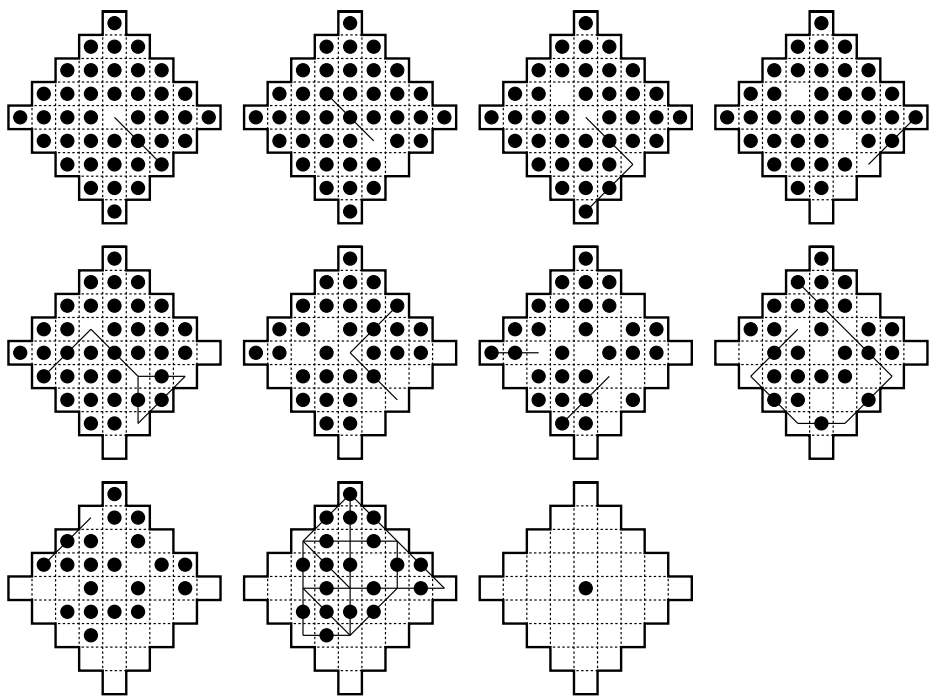}
\caption{A minimum length, 11-move solution to the central game on $Diamond(5)$.  This solution also contains the longest possible finishing move (18 jumps) for an 11-move solution to the central game.}
\label{fig11}
\end{figure}

\noindent\textbf{$Diamond(5)$ solutions}\newline
Central game in 11~moves (Figure~\ref{fig11}):
%g7-e5, d4-f6, i5-g7-e5, e9-g7, f4-h6-f8, c7-e9-g7, a5-c7, f2-d4-f6-d8-b6-d4, c3-a5-c5-c7-e5, h4-f4-d6-d8-f8-h6-h4-f2, e1-g3-e3-e1-c3-e3-e5-c3-c5-e5-g5-g7-e7-e5.
%Finishes with an 18-sweep
g7-e5, d4-f6, i5-g7-e5, e9-g7, b6-d4-f6-h6-f8-f6, g3-e5-g7, d8-f6, a5-c5, d2-f4-h6-f8-d8-b6-d4, b4-d2, e1-g3-i5-g5-e7-c7-c5-c3-e1-e3-g3-g5-e5-e7-c5-e5-c3-e3-e5.
$C9$~complement in 17~moves:
d2-f4, f8-d6, h4-f6, b6-d4, c3-e3-c5, g7-e5, e1-e3, h6-f6, b5-d5, e9-e7, i5-g5, a5-c3-c5, f7-f5, g4-e2-e4, c6-e8-e6, d8-b6-d4, f2-h4-f6.

\noindent\textbf{$Diamond(6)$ solution} (found using the techniques of Section 3.3)\newline
g6-complement in 15~moves (it is not known if this is the shortest possible):
i8-g6, f5-h7, i4-g6, k6-i4, h3-j5, g10-i8-k6-i4, f1-h3-j5, f7-h5-j7, d3-f5-h3, d9-f7, c6-e4-g2-i4-k6-i8-g6-e8-c6, b5-d3-f1-f3-h3-h5-j5-j7-h7-h9-f7-h7-h5-f5-f3-d3-d5-b5-d7-f7-f5, b7-d9-d7-d5-f7, f11-d9, a6-c6-c8-e10-g10-g8-e10-e8-g8-e6-e4-g6.

\vskip 30pt 
\footnotesize 
\begin{thebibliography}{mybib} 
\bibitem{Beasley} J. Beasley, \textit{The Ins and Outs of Peg Solitaire}, Oxford Univ. Press, paperback edition, 1992.
\bibitem{WinningWays} E. Berlekamp, J. Conway, and R. Guy, \textit{Winning Ways for your Mathematical Plays}, AK Peters, Volume 2: 695-734 (1982 edition) or Volume 4: 803--841 (2004 edition).
\bibitem{Bergholt} E. Bergholt, \textit{The Queen}, \textbf{131}: 666-667 (April 20) and 807 (May 11), 1912.
\bibitem{Bell} G. Bell, A fresh look at peg solitaire, \textit{Mathematics Magazine}, \textbf{80}:16--28, 2007.
\bibitem{GPJ28} J. Beasley (editor), \textit{The Games and Puzzles Journal}, Issue \textbf{28}, September 2003 \\http://www.gpj.connectfree.co.uk/gpjj.htm
\bibitem{BGS} G. Bell and J. Beasley, New problems on old solitaire boards,\newline http://arxiv.org/abs/math.CO/0611091, to appear in \textit{Board Game Studies}, Issue \textbf{8}, \newline http://www.boardgamesstudies.org/
\bibitem{Astar} P.E. Hart, N.J. Nilsson, and B. Raphael, A formal basis for heuristic determination of minimum path costs, \textit{IEEE Transactions on Systems Science and Cybernetics}, \textbf{4}:100--107, 1968.
\bibitem{Korf85} R. E. Korf, Depth-first iterative-deepening: An optimal admissible tree search, \textit{Artificial Intelligence}, \textbf{27}(1):97--109, 1985.
\bibitem{Korf96} R. E. Korf and L. A. Taylor, Finding optimal solutions to the twenty-four puzzle, in \textit{Proc. of the Int. Conf. on Artificial Intelligence (AAAI 96)}, 1202--1207, Portland, OR, 1996. 
\end{thebibliography}
\end{document}